# Comment on "On Nomenclature, and the Relative Merits of Two Formulations of Skew Distributions" by A. Azzalini, R. Browne, M. Genton, and P. McNicholas

Geoffrey J. McLachlan, Sharon X. Lee
Department of Mathematics, University of Queensland,
St. Lucia, Brisbane, Australia

**Abstract**

We comment on the recent paper by Azzalini et al. (2015) on two different distributions proposed in the literature for the modelling of data that have asymmetric and possibly long-tailed clusters. They are referred to as the restricted and unrestricted skew normal and skew $t$-distributions by Lee and McLachlan (2013a). We clarify an apparent misunderstanding in Azzalini et al. (2015) of this nomenclature to distinguish between these two models. Also, we note that McLachlan and Lee (2014) have obtained improved results for the unrestricted model over those reported in Azzalini et al. (2015) for the two datasets that were analysed by them to form the basis of their claims on the relative superiority of the restricted and unrestricted models. On this matter of the relative superiority of these two models, Lee and McLachlan (2014b, 2016) have shown how a distribution belonging to the broader class, the canonical fundamental skew $t$ (CFUST) class, can be fitted with little additional computational effort than for the unrestricted distribution. The CFUST class includes the restricted and unrestricted distributions as special cases. Thus the user now has the option of letting the data decide as to which model is appropriate for their particular dataset.

## 1 Introduction

In this paper, we provide some comments on Azzalini et al. (2015), which we shall refer to as ABGM in the sequel. In ABGM a comparison is given of two different distributions proposed for the modelling of data that have asymmetric and possibly long-tailed clusters. They refer to the two models as the classical and SDB, the latter so named since it was proposed by Sahu, Dey, and Branco (2003). These two distributions were referred to as the restricted multivariate skew $t$ (rMST) and unrestricted multivariate skew $t$ (uMST) distributions by Lee and McLachlan (2013a). We shall continue to use this latter terminology in our comments below.

In our comments we first wish to respond to statements in ABGM that are apparently based on a serious misunderstanding of the reporting of our results in Lee and McLachlan (2014a) and, in particular, of the nomenclature used therein. The discussion of our work in ABGM is limited to Lee and McLachlan (2014a), and so it does not consider the results presented in our other papers, in particular, Lee and McLachlan (2013a, 2013b, 2013c, 2013d, 2014b, 2015, 2016), although the first two of the latter seven papers are cited in ABGM. It is particularly unfortunate that these papers are not included in the comparison in ABGM as they contain a comparison of the restricted and unrestricted models applied to nine datasets from various fields.

Also, explicit cautionary notes are made in them to guard against any potential misunderstanding of our terminology. For example, in Lee and McLachlan (2013a, Page 244) it is stated that "Note that the use of 'restricted' here refers to restrictions on the random vector in the (conditioning-type) stochastic definition of the skew distribution. It is not a restriction on the parameter space, and so a 'restricted' form of a skew distribution is not necessarily nested within its corresponding 'unrestricted' form." In Lee and McLachlan (2013c), it is cautioned in the last two lines of Page 431 that "It should be stressed that the rMST family and uMST family match only in the univariate case, and one cannot obtain (7) [the restricted distribution] from (9) [the unrestricted distribution] when $p > 1$."

We also note that McLachlan and Lee (2014) have obtained improved results for the unrestricted model over those reported in Azzalini et al. (2015) and in their earlier paper Azzalini et al. (2014) for the two real datasets that were analysed by them to form the basis of their claims on the relative superiority of the restricted and unrestricted models.

The deficiencies in these two models have been demonstrated in Lee and McLachlan (2014b, 2016). Briefly, the restricted distribution is limited essentially to modelling skewness concentrated in a single direction in the feature space. This is because it is uses a univariate skewing function; that is, a single latent skewing variable is used in its convolution formulation. As a consequence, the realizations of the latent term used in the formulation of the model to represent skewness are confined to lie on a line in the $p$-dimensional feature space regardless of the value of $p$. This effectively means that the restricted distribution is limited to modelling skewness that is concentrated in a single direction in the feature space.

The unrestricted distribution on the other hand uses a multivariate skewing function with the feature-specific skewing variables that allow for skewness in the model taken to be uncorrelated. In its formulation, the $p$-dimensional vector of these skewing variables is premultiplied by a (diagonal) matrix of skewness parameters. The consequent net effect is that the feature-specific latent terms representing skewness in the model are uncorrelated. Thus it is suitable to model skew data with feature variables that are uncorrelated or approximately so. However, as the vector of skewing variables is premultiplied by a diagonal rather than an arbitrary matrix of skewness parameters, it may not model adequately skew data with strongly correlated features.

Lee and McLachlan (2014b, 2016) have shown how a distribution belonging to the broader class, the CFUST class, can be fitted with essentially no additional computational effort than for the unrestricted distribution. The CFUST distribution, which includes the restricted and unrestricted distributions as special cases, can model skewness in multivariate data with correlated feature variables since it uses an arbitrary matrix of skew parameters in its formulation.

With the availability of software for the fitting of mixtures of CFUST distributions (Lee and McLachlan, 2015), users now have the option for letting the data decide as to which model is appropriate for their particular dataset. Or they can fit all three models rMST, uMST, and CFUST (or mixtures of them), and make their own choice between the three using, say, an information-based criterion such as BIC.

Several of the claims in Azzalini et al. (2015) were made in their earlier paper Azzalini et al. (2014) and were commented on in McLachlan and Lee (2014).

## 2 Explanation of nomenclature for skew $t$-distributions

In an attempt to provide an automated approach to the clustering of flow cytometry data,

Pyne et al. (2009) considered the fitting of mixtures of skew $t$-distributions that belonged to the family of skew $t$-distributions proposed by Sahu et al. (2003). Members of the latter family have the following convolution-type characterization. The $p \times 1$ random vector $\boldsymbol{Y}$ can be expressed as

$$\boldsymbol{Y} = \boldsymbol{\mu} + \boldsymbol{\Delta}|\boldsymbol{U}_0| + \boldsymbol{U}_1, \qquad (1)$$

where

$$\begin{bmatrix} \boldsymbol{U}_0 \\ \boldsymbol{U}_1 \end{bmatrix} \sim N_{2p}\left(\begin{bmatrix} \boldsymbol{0} \\ \boldsymbol{0} \end{bmatrix}, \frac{1}{w}\begin{bmatrix} \boldsymbol{I}_p & \boldsymbol{0} \\ \boldsymbol{0} & \boldsymbol{\Sigma} \end{bmatrix}\right). \qquad (2)$$

In the above, $\boldsymbol{\mu}$ is a $p$-dimensional vector, $\boldsymbol{\Delta}$ is a $p \times p$ diagonal matrix, $\boldsymbol{I}_p$ denotes the $p \times p$ identity matrix, $\boldsymbol{\Sigma}$ is a positive definite matrix, and $\boldsymbol{0}$ is a vector/matrix of zeros with appropriate dimensions. Also, $w$ is the realization of the random variable $W$ distributed as gamma$(\frac{\nu}{2}, \frac{\nu}{2})$, and $|\boldsymbol{U}_0|$ denotes the vector whose $i$th element is the magnitude of the $i$th element of the vector $\boldsymbol{U}_0$.

In order to simplify the application of the EM algorithm to fit mixtures of these skew $t$-distributions, Pyne et al. (2009) imposed the restriction

$$U_{01} = U_{02} = \ldots = U_{0p} \qquad (3)$$

on the $p$ latent skewing variables, where $U_{0i} = (\boldsymbol{U}_0)_i$ ($i = 1, \ldots, p$). This produces a distribution equivalent to the skew $t$-distribution formulated by Branco and Dey (2001) and Azzalini and Capitanio (2003) after reparameterization. Lee and McLachlan (2013a) termed this distribution the restricted multivariate skew $t$ (rMST) distribution to distinguish it from the distribution proposed by Sahu et al. (2003). By default, the latter was referred to as the unrestricted multivariate skew $t$ (uMST) distribution since it can be characterized without any restrictions on the $p$ latent skewing variables in the convolution-type stochastic formulation (1).

By letting the degrees of freedom $\nu$ go to infinity in (1), we obtain a similar formulation for the restricted multivariate skew normal (rMSN) and unrestricted multivariate skew normal (uMSN) distributions. In the sequel we focus only on skew $t$-distributions since the situation is similar for skew normal distributions. One slight difference is that although the joint distribution of independent univariate skew normal random variables is the unrestricted skew normal distribution, the joint distribution of independent univariate skew $t$-random variables is only equal to the unrestricted skew $t$-distribution in the limit as the degrees of freedom $\nu$ in the marginal skew $t$-distributions becomes infinite.

## 3 Canonical fundamental skew $t$-distribution

The CFUST distribution was introduced (and so named) as a canonical version (special case) of the fundamental skew $t$-distribution by Arellano-Valle and Genton (2005). Its density is given by (1) where the vector $|\boldsymbol{U}_0|$ of latent skewing variables is taken to be of dimension $q(q \leq p)$ and where now the matrix $\boldsymbol{\Delta}$ of skewness parameters is a $p \times q$ matrix; $q$ is not necessarily restricted to being less than $p$.

An attractive feature of the CFUST distribution is that it includes the restricted and unrestricted distributions as special cases. And it can be fitted with little additional effort over the fitting of the unrestricted skew $t$-distributions (Lee and McLachlan, 2015, 2016).

We would like to point out that in our preprint Lee and McLachlan (2014b) and our subsequent paper Lee and McLachlan (2016), which was published online in February of last year, we have provided the EM equations for the fitting of a mixture of CFUST distributions. Lee and McLachlan (2015) have since given an R package for the fitting of a mixture of canonical fundamental skew $t$-distributions.

## 4 Unrestricted versus restricted skew $t$-distribution

The imposition of the restriction (3) on the $p$ latent skewing variables led Lee and McLachlan (2014a) to state that "the uMST distribution can be viewed as a simple extension of the rMST distribution." The use of "simple" there referred to the formulation of the models and was not meant to convey that the actual fitting of the unrestricted model is a relatively simple matter compared to the restricted model. The additional computational effort in fitting a mixture of uMST component distributions is perhaps best illustrated by the need for the calculation of the complicated multi-dimensional integrals arising on the E-step of the EM algorithm as described in Lee and McLachlan (2014b, 2016).

Also, we wish to stress that the use of "extension" in Lee and McLachlan (2014a) does not necessarily imply that the restricted distribution is a special case of the unrestricted version. Thus we do not follow the statement in ABGM that "Furthermore, the use of 'extension' is clearly inappropriate because neither one of the two families is a subset of the other for $p > 1$."

Concerning the theoretical results in ABGM to show that the amount of skewness in the rMST distribution is unlimited, we are not saying there is a limit on the amount of skewness under the restricted model; rather the limitation refers to having only a univariate skewing function. As a consequence, it is limited to modelling skewness concentrated in a single direction in the feature space.

More specifically, under the restriction (3) with the restricted skew $t$-distribution, realizations of the $p$-dimensional skewness term $\mathbf{\Delta}|\mathbf{U}_0|$ in the formulation (1) are confined to lie on a line regardless of the value of $p$. For example, in the case of $p = 2$ feature variables, the realizations of $\mathbf{\Delta}|\mathbf{U}_0|$ under the restriction (3) will lie on the line $y_2 = (\delta_2/\delta_1)y_1$, where $\delta_i = (\mathbf{\Delta})_{ii}$. It can be seen that the restriction (3) ensuring a univariate skewing function effectively means that the restricted $t$-distribution is limited to modelling skewness that is concentrated in a single direction in the feature space.

The unrestricted model can allow for skewness in more than one direction such as occurring with independent or uncorrelated skew feature variables. Although the unrestricted model has only the same number of skewness parameters as the restricted model, it has a multivariate skewing function with a $p$-dimensional latent vector of skewing variables in its formulation (1). Concerning the role of this skewing vector, it is claimed in ABGM on Page 5 that "In reality, the use of a multivariate latent error term in place of a single random component does not add any level of generality because this multivariate latent variable ... has a highly restricted structure." But this latent skewing vector ($|\mathbf{U}_0|$ in our notation here) does not have a highly restricted structure since $\mathbf{U}_0$ represents white noise (having mean zero and scale matrix equal to the identity matrix). However, as the absolute values of $\mathbf{U}_0$ are premultiplied by a diagonal matrix $\mathbf{\Delta}$ of skewness parameters in the formulation (1) of the unrestricted $t$-distribution, the feature-specific terms allowing for skewness in the feature variables are uncorrelated. Thus it is best suited to modelling skew data with uncorrelated features. For correlated feature data, a non-diagonal matrix $\mathbf{\Delta}$ is usually needed as illustrated in a series of examples in Lee and McLachlan (2014b, 2016) on the CFUST distribution.

# 5 Existence of improved fits for unrestricted skew $t$-mixtures

As the restricted skew $t$-distribution is not nested within the unrestricted skew $t$-family, one can generate datasets where one will be preferable to the other. To this end, Lee and McLachlan (2014b, 2016) have provided a series of examples to demonstrate situations where either one or both of the restricted and unrestricted skew $t$-distributions do not provide adequate models for skew data.

ABGM compared the relative clustering performance of the restricted and unrestricted skew $t$-distributions by using them to cluster two real datasets, the so-called crab dataset and the AIS dataset. It represented their claim that "The goal of the analyses herein is to present an extensive comparison ..." This would not appear to be an "extensive comparison" given only two datasets were considered and then only two- and three-dimensional subsets of them. Although the consideration of all possible such subsets resulted in 480 clusterings, they relate to only two datasets. Also, as explained in Section 4, the restricted model should be appropriate at least for bivariate skew data provided the features are not independent or weakly correlated.

Concerning the crab dataset as considered in ABGM, the correlations between any two of the five variables is very high (the lowest is 0.89). So this represents an ideal situation for the restricted model. More precisely, in the context of bivariate combinations of the variables for these data, it effectively means that bivariate datasets will lie almost on a straight line and so the restricted model with its univariate skewing function should not be disadvantaged. Thus on bivariate datasets of the crab dataset, the restricted model cannot perform too far below that of the unrestricted model and indeed should be better for those at least with highly correlated features.

In ABGM (Page 7), it is stated with respect to both the crab and AIS datasets that "The results (Figure 1) very clearly indicate that neither formulation is markedly superior and, if these results were to be taken in favour of either formulation, it would be the classical (restricted) formulation." However, in their analysis of all sets of the crab dataset, McLachlan and Lee (2014) found that the unrestricted version performs slightly better than the restricted. For example, on considering all 26 sets corresponding to the 26 different combinations of the five variables, they found that the restricted model gave a better fit for only 3 of the 26 sets. The differences were generally small with there being 13 ties. The EMMIXskew and EMMIXuskew packages were modified so that the two models could be started using the same values. The known class labels for these two datasets were not used as starting values in our cluster analyses.

As for the AIS dataset, it is not surprising that the performances of the restricted and unrestricted models are quite similar for the two- and three-dimensional combinations of the variables considered in ABGM, particularly for the bivariate combinations of the variables that have high correlations. On considering all 220 datasets corresponding to the 220 combinations consisting of all pairs and triplets of the 11 variables, McLachlan and Lee (2014) found that the unrestricted model gave a better fit for 105 versus 100 combinations for the restricted. The differences were generally small with there being 15 ties. But it is in contrast to the result above as reported in ABGM. McLachlan and Lee (2014) also found on using all of the 11 available variables that the unrestricted model gave a misclassification rate of 0.0198 compared to 0.0297 for the restricted model (that is, two fewer misallocations).

# 6 Concluding remarks

We have explained the motivation and the context in which Lee and McLachlan (2013a) adopted the terminology of restricted and unrestricted to describe two particular skew $t$-distributions. On their applicability as models, we have pointed out how Lee and McLachlan (2014b, 2016) have presented a series of examples to demonstrate the differences between the three models, the restricted, the unrestricted, and the so-called CFUST models.

Under the formulation (1) of these models, the restricted model is severely limited by having only a single latent skewing variable under its restriction (3). It is effectively limited to modelling skewness concentrated in a single direction in the feature space, whereas the unrestricted model allows for skewness in more than one direction such as with independent or uncorrelated feature variables.

The restricted and unrestricted models are special cases of the more general CFUST model, which also can model adequately multivariate data with correlated features since in its formulation it uses an arbitrary matrix of skewness parameters in conjunction with a multivariate skewing function.

We have also noted that McLachlan and Lee (2014) reported improved clustering performance of the unrestricted $t$-mixture model compared to that reported in ABGM for the crab dataset, which was one of two real datasets that they analysed. Concerning the other dataset analysed in ABGM, the AIS dataset, we note that McLachlan and Lee (2014) found for the two- and three-dimensional combinations of the variables in this dataset considered in ABGM that the unrestricted and restricted models perform very similarly but with the former slightly shading the restricted. This is in contrast to the result reported in ABGM.